\documentclass[12pt]{article}
\usepackage{amsmath}
\usepackage{mathrsfs}
\usepackage{bbm}
\usepackage{amsfonts,amsthm}
\usepackage{amssymb,a4}
\usepackage{extarrows}
\usepackage{color}
\usepackage[all,pdf]{xy}

	\allowdisplaybreaks  

\newcommand{\Z}{{\mathbb Z}}

\newcommand{\C}{{\mathbb C}}
\newcommand{\G}{\mathcal G}
\newcommand{\A}{\mathcal A}

\newcommand{\be}{\beta}
\newcommand{\dt}{\delta}
\newcommand{\lmd}{\lambda}
\newcommand{\gm}{\gamma}
\newcommand{\vf}{\varphi}
\newcommand{\e}{{\epsilon}}
\newcommand{\sgm}{\sigma}

\newcommand{\ing}{\in G}
\newcommand{\inc}{\in \C}
\newcommand{\der}{\text{Der}}
\newcommand{\ad}{\text{ad}}
\newcommand{\ov}{\overline}

\newcommand{\g}{\mathfrak g}

\newcommand{\setgl}[1]{G'(L_{#1})}
\newcommand{\setgi}[1]{G'(I_{#1})}
\newcommand{\cli}[1]{C_{LI}^{(#1)}}

\newcommand{\pf}[1]{\noindent{\bf Proof.}#1\hfill{}$\Box$}
\newcommand{\spanc}[1]{\text{span}_{\C}\{#1\}}

\newtheorem{thm}{Theorem}[section]
\newtheorem{prop}[thm]{Proposition}
\newtheorem{lem}[thm]{Lemma}

\begin{document}

\begin{center}
{\Large {\bf Deformed higher rank Heisenberg-Virasoro algebras}}\\
\vspace{0.5cm}
\end{center}

\begin{center}
{Chengkang Xu\footnote{
The author is supported by the National Natural Science Foundation of China (No. 11801375)}\\
Shangrao Normal University, Shangrao, Jiangxi, China\\
Email: xiaoxiongxu@126.com}
\end{center}

\begin{abstract}
In this paper, we study a class of infinitesimal deformations of the centerless higher rank Heisenberg-Virasoro algebras.
Explicitly, the universal central extensions, derivations and isomorphism classes of these algebras are determined.
\end{abstract}

\section{Introduction}

\def\theequation{1.\arabic{equation}}
\setcounter{equation}{0}

Let $G$ be an additive subgroup of $\C$ and $\lmd\inc$.
The Lie algebra $\g(G,\lmd)$ we study in this paper has a basis $\{L_a,I_a\mid a\ing\}$ subjecting to the following Lie brackets
\begin{equation}\label{eq1.1}
 [L_a,L_b]=(b-a)L_{a+b};\ \ \ [L_a,I_b]=(b-\lmd a)I_{a+b};\ \ \ \ [I_a,I_b]=0.
\end{equation}
We will simply denote $\g=\g(G,\lmd)$ if no confusion.
Since for any nonzero $\e\inc$, there is a Lie algebra isomorphism between $\g(G,\lmd)$ and $\g(\e G,\lmd)$ defined by
$$L_a\mapsto\e^{-1}L_{\e a};\ \ I_a\mapsto \e^{-1}I_{\e a},$$
we always assume in this paper that $\Z\subseteq G$ and $\frac1k\notin G$ for any integer $k>1$.

The algebra $\g$ is formed through semi-product of the higher rank Virasoro algebra $\mathfrak v=\spanc{L_a\mid a\ing}$
and one of its module of tensor fields $V(\lmd)=\spanc{v_a\mid a\ing}$ defined by
$$L_av_b=(b-\lmd a)v_{a+b}.$$
Notice that the algebra $\g(G,0)$ is the centerless higher rank Heisenberg-Virasoro algebra \cite{LZ},
and $\g(G,\lmd)$ is an infinitesimal deformation of $\g(G,0)$.
Therefore, we call $\g(G,\lmd)$ a {\em centerless deformed higher rank Heisenberg-Virasoro algebra}.
It also has close connections with the Heisenberg-Virasoro algebra and deformed Heisenberg-Virasoro algebras.
Taking $G=\Z$, we get the centerless deformed Heisenberg-Virasoro algebra $\g(\Z,\lmd)$,
which was given in \cite{LP} through the second cohomology group
of the centerless Heisenberg-Virasoro algebra $\g(\Z,\lmd)$
with coefficients in the adjoint representation.
Hence the deformed Heisenberg-Virasoro algebras and higher rank Heisenberg-Virasoro algebras are
generalizations of the Heisenberg-Virasoro algebra,
which was first introduced in \cite{ACKP}.
Representations of these three algebras were studied in \cite{L,LG,LvZ}
and references therein.

Due to the connection with these algebras,
it seems reasonable and important to study the algebra $\g$.
Our main purpose of this paper is to determine the universal central extensions,
derivations and isomorphism classes of all $\g(G,\lmd)$.
Such computation is a priority in the structure theory of Lie algebra,
and has been done for many Lie algebras,
such as generalized Witt algebras \cite{DZ},
generalized Schrodinger-Virasoro algebras \cite{WLX, TZ},
higher rank Heisenberg-Virasoro algebras \cite{LZ},
deformed Heisenberg-Virasoro algebras \cite{LP},
solenoidal Lie algebras over a quantum torus \cite{Xu},
and so on.

The following is the outline of this paper.
In section 2, we compute the universal central extension $\ov\g$ of $\g$,
which we call {\em deformed higher rank Heisenberg-Virasoro algebra}.
Section 3 is devoted to derivations of $\g$ and the lifts of these derivations
to derivations of $\ov\g$.
In the last section isomorphism classes and automorphism groups of
all the algebras $\g(G,\lmd)$ are determined,
and finally the explicit form of the unique lift of any outer automorphism of $\g$ to $\ov\g$
is given.

Throughout this paper, the symbols $\Z,\C,\C^\times$ refer to the set of integers,
complex numbers and nonzero complex numbers respectively.
We denote by $G^*$ the set of characters of $G$,
and by $\text{Hom}(G,\C)$ the set of all additive group homomorphisms.

\section{Universal central extensions}\label{sec2}
\def\theequation{2.\arabic{equation}}
\setcounter{equation}{0}

In this section we compute the universal central extensions of $\g$
for $\lmd\neq-1$, $G$ being free and of rank $n\geq1$.
We fix a $\Z$-basis $\e_1=1,\e_2,\cdots,\e_n$ of $G$.

Let $\vf:\g\times\g\longrightarrow\C$ be a 2-cocycle, hence
\begin{equation}\label{eq2.01}
 \begin{split}
  &\vf([x,y],z)+\vf([y,z],x)+\vf([z,x],y)=0;\\
  &\vf(x,y)=-\vf(y,x),\ \ \text{ for any }x,y,z\in\g.
 \end{split}
\end{equation}
Define a linear function $f$ on $\g$ by
$$\begin{aligned}
 &f(L_a)=\frac1 a\vf(L_0,L_a);\ \ f(I_a) =\frac1 a\vf(L_0,I_a) \text{ for }a\neq0;\\
 &f(L_0)=\frac12\vf(L_{-1},L_1);\ \ f(I_0)=\frac1{\lmd+1}\vf(L_{-1},I_1),
\end{aligned}$$
and a 2-coboundary $\vf_f$ by $\vf_f(x,y)=f([x,y])$.
Replacing $\vf$ by $\vf-\vf_f$, we may assume
$$\vf(L_0,L_a)=\vf(L_0,I_a)=\vf(L_{-1},L_1)=\vf(L_{-1},I_1)=0\ \ \text{ for any }a\neq0.$$

\begin{lem}\label{lem2.1}
 $\vf(L_a,L_b)=\frac1{12}(a^3-a)C_L\dt_{a+b,0}$ for some $C_L\in\C$.
\end{lem}
\pf{
One reference for this computation is \cite{PZ}.
}

\begin{lem}\label{lem2.2}
 $\vf(I_a,I_b)=aC_I\dt_{a+b,0}\dt_{\lmd,0}$ for some $C_I\in\C$.
\end{lem}
\pf{
Consider the triple $(x,y,z)=(L_a,I_b,I_c)$ in (\ref{eq2.01}) and we have
\begin{equation}\label{eq2.02}
 (b-\lmd a)\vf(I_{a+b},I_c)+(c-\lmd a)\vf(I_b,I_{a+c})=0.
\end{equation}
Set $a=0$ in (\ref{eq2.02}) and we get
$$\vf(I_b,I_c)=0 \text{ if }b+c\neq0.$$
Take $c=-a-b$ in (\ref{eq2.02}),
\begin{equation}\label{eq2.03}
 (b-\lmd a)\vf(I_{a+b},I_{-a-b})=((\lmd+1)a+b)\vf(I_b,I_{-b}).
\end{equation}
Set $a=1-b$ in (\ref{eq2.03}), we obtain
\begin{equation}\label{eq2.04}
 -(\lmd b-\lmd-1)\vf(I_b,I_{-b})=((\lmd+1)b-\lmd)\vf(I_1,I_{-1}),
\end{equation}
which implies the lemma for $\lmd=0$ by setting $C_I=\vf(I_1,I_{-1})$.

Now we only need to prove that $\vf(I_a,I_{-a})=0$ for any $a\ing$ if $\lmd\neq0,-1$.
If $\lmd^{-1}\in\Z$, then set $a=\frac b \lmd \ing$ in (\ref{eq2.03}) and one gets
$$b(\frac1\lmd+2)\vf(I_b,I_{-b})=0,$$
which shows that $\vf(I_b,I_{-b})=0$ if $\lmd\neq-\frac12$.
Suppose $\lmd=-\frac12$ and by (\ref{eq2.04}) we see
$$\vf(I_b,I_{-b})=\vf(I_1,I_{-1}) \text{ if }b\neq-1.$$
Choose $b\neq\pm1$ and let $a=-1-b$ in (\ref{eq2.03}), we get
$$\frac{b-1}2\vf(I_{-1},I_1)=\frac{b-1}2\vf(I_b,I_{-b})=\frac{b-1}2\vf(I_1,I_{-1}),$$
which forces that $\vf(I_1,I_{-1})=0$ by skew-symmetry.
So $\vf(I_b,I_{-b})=0$.

Suppose that $\lmd^{-1}\notin\Z$, then $\lmd b+\lmd+1\neq0$ for any $b\ing$.
So (\ref{eq2.04}) gives
$$\vf(I_a,I_{-a})=-\frac{(\lmd+1)a-\lmd}{\lmd a-\lmd-1}\vf(I_1,I_{-1}).$$
Put it into (\ref{eq2.03}) and set $b=0$, then we get
$$\lmd(\lmd+1)(2\lmd+1)a(a-1)\vf(I_1,I_{-1})=0\text{ for any }a\ing,$$
which forces $\vf(I_1,I_{-1})=0$, and hence $\vf(I_a,I_{-a})=0$.
We complete the proof.
}

\begin{lem}\label{lem2.3}
 $\vf(L_a,I_b)=\dt_{a+b,0}\left(C_{LI}^{(0)}(a^2+a)\dt_{\lmd,0}
  +\frac1{12}(a^3-a)C_{LI}^{(1)}\dt_{\lmd,1}
  +\sum\limits_{i=2}^na_iC_{LI}^{(i)}\dt_{\lmd,-2}\right)$,
  where $C_{LI}^{(i)}\in\C, 0\leq i\leq n$ and $a=\sum\limits_{j=1}^na_j\e_j$.
\end{lem}
\pf{
Consider the triple $(x,y,z)=(L_a,L_b,I_c)$ in (\ref{eq2.01}), we have
\begin{equation}\label{eq2.05}
 (b-a)\vf(L_{a+b},I_c)-(c-\lmd b)\vf(L_a,I_{b+c})+(c-\lmd a)\vf(L_b,I_{a+c})=0.
\end{equation}
Take $a=0$,
$$(b+c)\vf(L_b,I_c)-(c-\lmd b)\vf(L_0,I_{b+c})=0,$$
which implies
\begin{align}\label{eq2.06}
\vf(L_b,I_c)=0\ \text{ if }\ b+c\neq0.
\end{align}
Moreover, let $b+c=0$ and $c-\lmd b\neq0$ (such $b,c$ exist since $\lmd\neq-1$),
we get $\vf(L_0,I_0)=0$.

Set $c=-a-b$ in (\ref{eq2.05}) and we obtain
\begin{equation}\label{eq2.07}
 (b-a)\vf(L_{a+b},I_{-a-b})+(a+b(\lmd+1))\vf(L_a,I_{-a})-(a(\lmd+1)+b)\vf(L_b,I_{-b})=0.
\end{equation}
Let $a+b=0$ in (\ref{eq2.07}) and we get
\begin{equation}\label{eq2.08}
 \vf(L_{-a},I_a)=-\vf(L_a,I_{-a})\text{ if }\lmd\neq0.
\end{equation}
Let $b=1$ in (\ref{eq2.07}), we have
\begin{equation}\label{eq2.09}
 (1-a)\vf(L_{a+1},I_{-a-1})+(a+\lmd+1)\vf(L_a,I_{-a})-(a(\lmd+1)+1)\vf(L_1,I_{-1})=0.
\end{equation}
Moreover, setting $b=-1$ and replacing $a$ by $a+1$ in (\ref{eq2.07}), it gives
\begin{equation}\label{eq2.10}
 (a-\lmd)\vf(L_{a+1},I_{-a-1})-(a+2)\vf(L_a,I_{-a})=0.
\end{equation}
Now we continue the proof in the following four cases.

{\bf Case I}: $\lmd=0$.
Combining (\ref{eq2.09}) and (\ref{eq2.10}), we get
$$\vf(L_a,I_{-a})=\frac12(a^2+a)\vf(L_1,I_{-1}).$$
Set $C_{LI}^{(0)}=\frac12\vf(L_1,I_{-1})$,
and we get the lemma by (\ref{eq2.06}).

We emphasize that in the remaining three cases we have
$\vf(L_{-a},I_a)=-\vf(L_a,I_{-a})$ for any $a\ing$.
Especially, $\vf(L_1,I_{-1})=-\vf(L_{-1},I_1)=0$.

{\bf Case II}: $\lmd\neq 0,\pm1,-2$.
Notice that (\ref{eq2.09}) and (\ref{eq2.10}) form a system of homogeneous liner equations
in variables $\vf(L_a,I_{-a})$ and $\vf(L_{a+1},I_{-a-1})$ with a coefficient matrix
whose determinant is
$$\left|
  \begin{array}{cc}
    1-a    & a+\lmd+1 \\
    a-\lmd & -a-2\\
  \end{array}
\right|=(\lmd-1)(\lmd+2)\neq0.$$
So $\vf(L_a,I_{-a})=0$ and hence by (\ref{eq2.06}) $\vf(L_a,I_b)=0$ for any $a,b\ing$.

{\bf Case III}: $\lmd=1$.
Set $C_{LI}^{(1)}=2\vf(L_2,I_{-2})$.
A same calculation as in Lemma \ref{lem2.1} shows
$$\vf(L_a,I_b)=\frac1{12}(a^3-a)C_{LI}^{(1)}\dt_{a+b,0}.$$

{\bf Case IV}: $\lmd=-2$.
In this case, the equations (\ref{eq2.09}) and (\ref{eq2.10}) turn to
$$\begin{aligned}
 (a-1)\vf(L_a,I_{-a})&=(a-1)\vf(L_{a+1},I_{-a-1}),\\
 (a+2)\vf(L_a,I_{-a})&=(a+2)\vf(L_{a+1},I_{-a-1}),
\end{aligned}$$
which implies that $\vf(L_{a+1},I_{-a-1})=\vf(L_a,I_{-a})$ for any $a\ing$. So
$$\vf(L_k,I_{-k})=0,\ \ \vf(L_{a+k},I_{-a-k})=\vf(L_a,I_{-a})
\text{ for any }k\in\Z, a\ing\setminus\Z.$$
Then (\ref{eq2.05}) gives
$$\vf(L_{a+b},I_{-a-b})=\vf(L_a,I_{-a})+\vf(L_b,I_{-b})\text{ if }a\neq b.$$
For $a\ing\setminus\Z$, we still have
$$\vf(L_{2a},I_{-2a})=\vf(L_{a+1},I_{-a-1})+\vf(L_{a-1},I_{-a+1})=2\vf(L_a,I_{-a}).$$
So we have proved $\vf(L_{a+b},I_{-a-b})=\vf(L_a,I_{-a})+\vf(L_b,I_{-b})$ for any $a,b\ing.$
Therefore $\vf(L_a,I_{-a})=\sum\limits_{i=2}^na_i\vf(L_{\e_i},I_{-\e_i})$
if $a=\sum\limits_{i=1}^na_i\e_i$.
Put $C_{LI}^{(i)}=\vf(L_{\e_i},I_{-\e_i})$ for $2\leq i\leq n$,
and we see that $\vf(L_a,I_b)=\sum\limits_{i=2}^na_iC_{LI}^{(i)}\dt_{a+b,0}$.
This finishes the proof.
}

Now from Lemma \ref{lem2.1}, \ref{lem2.2} and \ref{lem2.3} we get the main theorem of this section.
\begin{thm}\label{thm2.4}
Denote by $\overline\g$ the universal central extension of $\g(G,\lmd)$ for $\lmd\neq-1$
and $G$ be a free subgroup of $\C$ of rank $n$.
Then $\overline\g$ satisfies the following Lie brackets
$$\begin{aligned}
 &[L_a,L_b] =(b-a)L_{a+b}+\frac1{12}(a^3-a)C_L\dt_{a+b,0};\ \ \ \
  [I_a,I_b] =aC_I\dt_{a+b,0}\dt_{\lmd,0};\\
 &[L_a,I_b] =(b-\lmd a)I_{a+b}+\dt_{a+b,0}\left(C_{LI}^{(0)}(a^2+a)\dt_{\lmd,0}
  +\frac1{12}(a^3-a)C_{LI}^{(1)}\dt_{\lmd,1}
  +\sum\limits_{i=2}^na_iC_{LI}^{(i)}\dt_{\lmd,-2}\right),
\end{aligned}$$
where $C_L,C_I,C_{LI}^{(i)},0\leq i\leq n$ are central elements and $a=\sum\limits_{i=1}^na_i\e_i$.
\end{thm}

\noindent
{\bf Remark}: (1) If $\lmd\neq-1,-2$,
Theorem \ref{thm2.4} stands for arbitrary additive subgroup $G$ of $\C$.
Especially, when $\lmd=0$, we get the generalized Heisenberg-Virasoro algebras,
which was originally given in \cite{LZ}.
When $G=\Z$, we get the deformed Heisenberg-Virasoro algebra,
which was studied in \cite{LP}.\\
(2) The Lie algebra $\g=\g(G,-1)$ is not perfect ($I_0\notin[\g,\g]$)
and there is no universal central extension of $\g$.
But one can still consider the universal central extension of the derived subalgebra
$\g'=[\g,\g]=\spanc{L_a,I_b\mid a\ing,b\ing\setminus\{0\}}$.
One can prove that
the universal central extension of $\g'$ is governed by four nontrivial 2-cocycles
$$\begin{aligned}
 C_L(L_a,L_b)&=\frac1{12}(a^3-a)\dt_{a+b,0};\ \ \
                         C_I(I_a,I_b)=\frac1a\dt_{a+b,0};\\
 C_{LI}(L_a,I_b)&=a\dt_{a+b,0};\ \ \ C_{LI}'(L_a,I_b)=\dt_{a+b,0}.
\end{aligned}$$

\section{Derivations of $\g$ and $\overline\g$}\label{sec3}
\def\theequation{3.\arabic{equation}}
\setcounter{equation}{0}

In this section we compute derivations of $\g$ for arbitrary $G$ and $\lmd$,
and consider the lifts of these derivations of $\g$ to derivations of $\overline\g$.

The algebra $\g$ has a natural $G$-grading $\g=\bigoplus\limits_{a\ing}\g_a$,
where $\g_a=\spanc{L_a,I_a}$.
From a well known result about derivations of graded Lie algebras in \cite{Far}
we know that the algebra $\der\g$ of derivations of $\g$ is also graded by $G$,
i.e, $\der\g=\bigoplus\limits_{a\ing}(\der\g)_a$, where
$$(\der\g)_a=\spanc{\sgm\in\der\g\mid\sgm(\g_b)\subseteq\g_{a+b},\ \forall\ b\ing},$$
and all outer derivations of $\g$ lie in $(\der\g)_0$.
Now we construct some derivations of $\g$ which are all of degree 0.
Define linear maps $\vf, \psi:\g\longrightarrow\g$ by
$$\vf(L_a)=aI_a,\ \vf(I_a)=0;\hspace{3cm}
   \psi(L_a)=0,\ \psi(I_a)=I_a.$$
For $\lmd=0$, we define a linear map $\sgm_{(0)}:\g\longrightarrow\g$ by
$$\sgm_{(0)}(L_a)=I_a\dt_{\lmd,0},\ \sgm_{(0)}(I_a)=0.$$
For $\lmd=-1$, we define a linear map $\sgm_{(-1)}:\g\longrightarrow\g$ by
$$\sgm_{(-1)}(L_a)=a^2I_a\dt_{\lmd,-1},\ \sgm_{(-1)}(I_a)=0.$$
For $\lmd=-2$, we define a linear map $\sgm_{(-2)}:\g\longrightarrow\g$ by
$$\sgm_{(-2)}(L_a)=a^3I_a\dt_{\lmd,-2},\ \sgm_{(-2)}(I_a)=0.$$
Let $\A\in\text{Hom}(G,\C)$.
We define a linear map $\xi_\A:\g\longrightarrow\g$ by
$$\xi_\mathcal{A}(L_a)=\mathcal{A}(a)L_a,\ \xi_\mathcal{A}(I_a)=\mathcal{A}(a)I_a,$$
and for $\lmd=1$, a linear map $\eta_{\A,1}:\g\longrightarrow\g$ by
$$\eta_{\A,1}(L_a)=\mathcal{A}(a)I_a\dt_{\lmd,1},\ \eta_{\A,1}(I_a)=0.$$
It is easy to check that the linear maps
$\vf,\psi,\sgm_{(0)},\sgm_{(-1)},\sgm_{(-2)},\xi_\A,\eta_{\A,1}$ are all
derivations of degree 0.
The main result of this section is the following

\begin{thm}\label{thm3.2}
 For $a\neq0$, $(\der\g)_a=\spanc{\text{ad}L_a,\text{ad}I_a}$, and
 $$(\der\g)_0=\spanc{\vf,\psi,\sgm_{(0)},\sgm_{(-1)},\sgm_{(-2)},\xi_\A,\eta_{\A,1}
                     \mid \mathcal A\in\text{Hom}(G,\C)}.$$
\end{thm}
\pf{
Let $\sgm\in(\der\g)_0$ and suppose
$$\sgm(L_a)=\mu_L(a)L_a+\tau_L(a)I_a,\ \ \sgm(I_a)=\mu_I(a)L_a+\tau_I(a)I_a,$$
for some functions $\mu_L,\tau_L,\mu_I,\tau_I$ on $G$.
Apply $\sgm$ to $[L_a,L_b], [L_a,I_b]$ and $[I_a,I_b]$, we get the following equations
\begin{align}
 &\mu_I(a)(b-\lmd a)=\mu_I(b)(a-\lmd b),\label{eq3.01}\\
 &\mu_I(a+b)(b-\lmd a)=\mu_I(b)(b-a),\label{eq3.02}\\
 &(b-\lmd a)\tau_I(a+b)=(b-\lmd a)\mu_L(a)+(b-\lmd a)\tau_I(b),\label{eq3.03}\\
 &(b-a)\tau_L(a+b)=(b-\lmd a)\tau_L(b)-(a-\lmd b)\tau_L(a),\label{eq3.04}\\
 &(b-a)\mu_L(a+b)=(b-a)(\mu_L(a)+\mu_L(b)).\label{eq3.05}
\end{align}

{\bf Claim 1}: $\mu_L\in\text{Hom}(G,\C)$.
From (\ref{eq3.05}) we see that $\mu_L(a+b)=\mu_L(a)+\mu_L(b)$ if $a\neq b$.
Moreover, let $a=0, b\neq0$ in (\ref{eq3.05}) we get $\mu_L(0)=0$.
Hence $\mu_L(-a)=-\mu_L(a)$ for any $a\ing$.
Then choose $b\neq 0,\pm a$ and we have
$$\mu_L(2a)=\mu_L(a+b)+\mu_L(a-b)=2\mu_L(a)+\mu_L(-b)+\mu_L(b)=2\mu_L(a).$$
So we have proved $\mu_L(a+b)=\mu_L(a)+\mu_L(b)$ for any $a,b\ing$.
Claim 1 stands.

{\bf Claim 2}: $\mu_I=0$.
Let $b=0$ in (\ref{eq3.02}) we get
\begin{equation}\label{eq3.06}
 \lmd\mu_I(a)=\mu_I(0)\text{ for any }a\neq0.
\end{equation}
If $\lmd=0$, then $\mu_I(0)=0$.
Take $b=-a$ in (\ref{eq3.02}) and we have $\mu_I(a)=0$ for any $a\ing$.
Suppose $\lmd\neq0$, then $\mu_I(a)=\frac1\lmd\mu_I(0)$ for $a\neq 0$.
Then (\ref{eq3.01}) turns to
$$\mu_I(0)(b-a)(\lmd+1)=0,$$
which implies $\mu_I(0)=0$ and hence $\mu_I=0$ if $\lmd\neq -1$.
If $\lmd=-1$, then take $a=-b\neq 0$ in (\ref{eq3.02}) and we have
$\mu_I(b)=0$ for any $b\neq 0$.
So $\mu_I(0)=0$ by (\ref{eq3.06}).
This proves Claim 2.

{\bf Claim 3}: $\tau_I(a)=\mu_L(a)+\tau_I(0)$.
If $\lmd=0$, then (\ref{eq3.03}) shows that
$$\tau_I(a+b)=\mu_L(a)+\tau_I(b)\text{ for }b\neq 0.$$
Choose $a=-b\neq 0$ then we have
$$\tau_I(0)=\mu_L(-b)+\tau_I(b)=-\mu_L(b)+\tau_I(b),$$
that is, $\tau_I(b)=\mu_L(b)+\tau_I(0)$.
This proves the claim for $\lmd=0$ since $\mu_L(0)=0$.
Suppose $\lmd\neq 0$. Let $b=0,a\neq 0$ in (\ref{eq3.03}),
and Claim 3 follows.

{\bf Claim 4}: $\tau_L(a)=l_0\dt_{\lmd,0}+l_1a+l_2a^2\dt_{\lmd,-1}+l_3a^3\dt_{\lmd,-2}+
        \mathcal B(a)\dt_{\lmd,1}$ for some $l_0,l_1,l_2,l_3\inc$
        and $\mathcal B\in\text{Hom}(G,\C)$.
Take $a=0$ in (\ref{eq3.04}) and we get
\begin{equation}\label{eq3.07}
 \lmd\tau_L(0)=0.
\end{equation}
Replacing $b$ by $-b$ and $a$ by $a+b$ in (\ref{eq3.04}) we have
\begin{equation}\label{eq3.08}
 (a+b(\lmd+1))\tau_L(a+b)-(a+2b)\tau_L(a)=-(\lmd a+b(\lmd+1))\tau_L(-b).
\end{equation}
Combining (\ref{eq3.04}) and (\ref{eq3.08}), and letting $b=1$, we obtain
\begin{equation}\label{eq3.09}
 (\lmd-1)(\lmd+2)\tau_L(a)=(\lmd a-1)(a+\lmd+1)\tau_L(1)+(a-1)(\lmd a+\lmd+1)\tau_L(-1).
\end{equation}
If $\lmd=0$, then we get
$$\tau_L(a)=l_1 a+l_0\ \text{ with }\
l_1=\frac12(\tau_L(1)-\tau_L(-1)),l_0=\frac12(\tau_L(1)+\tau_L(-1)).$$
If $\lmd\neq 0$ then $\tau_L(0)=0$ by (\ref{eq3.07}).
Let $a+b=0$ in (\ref{eq3.04}) and we get
$$\tau_L(-a)=-\tau_L(a)\text{ if }\lmd\neq 0,-1.$$
So by (\ref{eq3.09}) we see $\tau_L(a)=a\tau_L(1)$ if $\lmd\neq 0,\pm1,-2$.

If $\lmd=-1$, then from (\ref{eq3.09}) it follows
$$\tau_L(a)=l_2 a^2+l_1a\ \text{ with }\
  l_2=\frac12(\tau_L(1)+\tau_L(-1)),l_1=\frac 12(\tau_L(1)-\tau_L(-1)).$$
If $\lmd=1$, then (\ref{eq3.04}) and (\ref{eq3.08}) imply that
$$\tau_L(a+b)=\tau_L(a)+\tau_L(b)\text{ for any }a,b\ing,$$
that is, $\tau_L \in\text{Hom}(G,\C)$.
We may write $\tau_L(a)=\mathcal B(a)+l_1a$ for some $l_1\in\C$ and $\mathcal B\in\text{Hom}(G,\C)$.

Now suppose $\lmd=-2$.
Let $b=1$ in (\ref{eq3.04}),
$$(1-a)\tau_L(a+1)+(a+2)\tau_L(a)=(2a+1)\tau_L(1).$$
Furthermore, replace $a$ by $a+1$ and we have
$$-a\tau_L(a+2)+(a+3)\tau_L(a+1)=(2a+3)\tau_L(1).$$
Let $b=2$ in (\ref{eq3.04}),
$$(a-2)\tau_L(a+2)-(a+4)\tau_L(a)=-2(a+1)\tau_L(2).$$
Then the above three equations imply that
$$\tau_L(a)=l_3 a^3+l_1a\ \text{ with }\
 l_3=\frac16(\tau_L(2)-2\tau_L(1)),l_1=\frac 16(8\tau_L(1)-\tau_L(2)).$$
Combing the above cases we see that Claim 4 is valid.

In conclusion, we have
$$\begin{aligned}
 &\sgm(L_a)=\mathcal A(a)L_a+\left(l_0\dt_{\lmd,0}+l_1a+l_2a^2\dt_{\lmd,-1}+l_3a^3\dt_{\lmd,-2}+
        \mathcal B(a)\dt_{\lmd,1}\right)I_a;\\
 &\sgm(I_a)=\left(\mathcal A(a)+l\right)I_a,
\end{aligned}$$
for some $l,l_0,l_1,l_2,l_3\in\C$ and $\mathcal{A,B}\in\text{Hom}(G,\C)$.
This proves Theorem \ref{thm3.2}.
}

\noindent
{\bf Remark}:
(1) Denote by $\mathrm{id}:G\longrightarrow\C$ be the identity map.
The corresponding derivation $\xi_{\mathrm{id}}$ is exactly the inner derivation $\text{ad}L_0$.
The inner derivation $\text{ad}I_0=0$ if $\lmd=0$,
and $\text{ad}I_0=\frac1\lmd\vf$ if $\lmd\neq0$.\\
(2) In \cite{LZ}, derivations for higher rank Heisenberg-Virasoro algebras were computed.
However, there is one derivation missed, the one we denoted by $\sgm_{(0)}$.
The reason why they missed this derivation is that they assumed $\tau_L(0)=0$ out of nowhere
(or in the terminology of \cite{LZ} $\be_0=0$, see Line 3 Page 9 in \cite{LZ}).

In the following we consider lifts of the derivations of $\g$ we obtained above
to derivations of $\overline\g$ with $\lmd\neq-1$ and $G$ being free of rank $n\geq1$.
We shall recall a result about derivations of the universal central extension
of a perfect Lie algebra from \cite{BM}.
\begin{prop}[\cite{BM}]\label{prop3.3}
 Suppose the Lie algebra $\G$ is perfect and denote by $\overline\G$
 the universal central extension of $\G$.
 Then every derivation of $\G$ lifts to a derivation of $\overline\G$.
 Moreover, if $\G$ is centerless, then this lift is unique and $\der\overline\G\cong\der\G$.
\end{prop}

Denote $\mathfrak c=\spanc{C_L,C_I\dt_{\lmd,0},\cli 0\dt_{\lmd,0},\cli 1\dt_{\lmd,1},
               \cli i\dt_{\lmd,-2}\mid 2\leq i\leq n}$.
From Proposition \ref{prop3.3} we know that if $\lmd\neq0,-1$,
every derivation of $\g$ lifts uniquely to a derivation of $\overline\g$.
For any $\sgm\in\der\g$,
define a linear map $\overline\sgm:\overline\g\longrightarrow\overline\g$ by
\begin{align}\label{eq3.10}
 \overline\sgm(L_a)=\sgm(L_a);\  \overline\sgm(I_a)=\sgm(I_a);\  \overline\sgm(\mathfrak c)=0\ \
 \text{ for any }\ a\ing.
\end{align}
Clearly, $\overline\sgm$ is a derivation of $\overline\g$ and it lifts $\sgm$.
Therefore we get
\begin{thm}
 If $\lmd\neq 0,-1$, then
 $$\der\ov\g=\spanc{\ov\vf,\ov\psi,\ov{\sgm_{(-2)}},\ov{\xi_{\mathcal A}},
  \ov{\eta_{\mathcal A,1}},\ad L_a,\ad I_a\mid \mathcal A\in\text{Hom}(G,\C),a\ing},$$
 where $\vf,\psi,\sgm_{(-2)},\xi_\A, \eta_{\A,1}$ are as in Theorem \ref{thm3.2} and
 $\ov\vf,\ov\psi,\ov{\sgm_{(-2)}},\ov{\xi_\A}, \ov{\eta_{\A,1}}$ defined in (\ref{eq3.10}).
\end{thm}

Now we consider the $\lmd=0$ case.
We emphasize again that $\ov\g$ is the higher rank Heisenberg-Virasoro algebra in this case.
There are four kinds of derivations $\vf,\psi,\sgm_{(0)},\xi_\A$,
which are of degree 0.

\begin{thm}
Let $\lmd=0$ and $\A\in\text{Hom}(G,\C)$ (here we denote $C_{LI}=\cli0$).\\
(1) The derivation $\vf$ lifts uniquely to a derivation $\ov\vf$ of $\ov\g$ defined by
    \begin{align*}
     L_a\mapsto aI_a+\dt_{a,0}C_{LI},\ I_a\mapsto C_I\dt_{a,0},\ C_L\mapsto-24C_{LI},\
      C_{LI}\mapsto C_I,\  C_I\mapsto0.
    \end{align*}
(2) The derivation $\sgm_{(0)}$ lifts uniquely to a derivation $\ov{\sgm_{(0)}}$ of $\ov\g$
    defined by
    \begin{align*}
     L_a\mapsto I_a-\dt_{a,0}C_{LI},\ I_a\mapsto -C_I\dt_{a,0},\ C_L\mapsto0,\
      C_{LI}\mapsto0,\  C_I\mapsto0.
    \end{align*}
(3) The derivation $\xi_\A$ lifts to a family of derivations $\{\ov{\xi_{\A,l,k}}\mid l,k\inc\}$
    of $\ov\g$ defined by
    \begin{align*}
     &L_a\mapsto\A(a)L_a+(l+ka)I_a+\dt_{a,0}(k-l)C_{LI},\\
     &I_a\mapsto\A(a)I_a-\dt_{a,0}C_I(l+k),\\
     &C_L\mapsto-24kC_{LI},\ C_{LI}\mapsto-kC_I,\  C_I\mapsto0.
    \end{align*}
(4) The derivation $\psi$ lifts to a family of derivations $\{\ov{\psi_{l,k}}\mid l,k\inc\}$
    of $\ov\g$ defined by
    \begin{align*}
     &L_a\mapsto (l+ka)I_a+\dt_{a,0}(k-l)C_{LI},\\
     &I_a\mapsto I_a+\dt_{a,0}C_I(k-l),\\
     &C_L\mapsto-24kC_{LI},\ C_{LI}\mapsto C_{LI}+kC_I,\  C_I\mapsto2C_I.
    \end{align*}
\end{thm}
\pf{
We only prove (3), the other three are similar and we omit it.

Let $\phi$ be a lift of $\xi_\A$. Notice that $\phi$ is homogeneous of degree 0.
We may write
$$\phi(L_a)=\A(a)L_a+f_L(a)I_a+\dt_{a,0}C_1;\ \ \ \
  \phi(I_a)=f_I(a)L_a+\A(a)I_a+\dt_{a,0}C_2,$$
for some $C_1,C_2\in\mathfrak c$ and functions $f_L,f_I:G\longrightarrow\C$.
Expanding the equation $\phi([L_a,L_b])=[\phi(L_a),L_b]+[L_a,\phi(L_b)]$
and comparing coefficients we get
$(b-a)f_L(a+b)=bf_L(b)-af_L(a)$, which implies
$$f_L(a)=l+ka \ \text{ for some }\ l,k\inc,$$
and moreover we get
$$\frac1{12}(a^3-a)\phi(C_L)-2aC_1 =\left(f_L(-a)(a^2+a)-f_L(a)(a^2-a)\right)C_{LI}
  =2a(l-ka^2)C_{LI}.$$
Let $a=1$, we see $C_1=(k-l)C_{LI}$.
Hence $\phi(C_L)=-24kC_{LI}$.

Expanding $\phi([L_a,I_b])=[\phi(L_a),I_b]+[L_a,\phi(I_b)]$ we get
$bf_I(a+b)=f_I(b)(b-a)$, which implies $f_I(a)=0$ for any $a\ing$,
and we get
$$(a^2+a)\phi(C_{LI})-aC_2=a(l-k)C_I.$$
Let $a=-1$ and we see that $C_2=-(l+k)C_I,\phi(C_{LI})=-kC_I.$

At last expand $\phi([I_a,I_b])$ and we have $\phi(C_I)=0$.
This proves (3).
}

\section{Isomorphism classes and automorphisms}\label{sec4}
\def\theequation{4.\arabic{equation}}
\setcounter{equation}{0}

In this section we determine the isomorphism classes
and automorphism groups of the Lie algebras $\g(G,\lmd)$ for arbitrary $G$ and $\lmd$,
and then consider the lifts of these automorphisms to
automorphisms of $\ov\g$ for $\lmd\neq-1$ and free $G$.
\begin{thm}\label{thm4.1}
 The Lie algebras $\g(G,\lmd)$ and $\g(G',\lmd')$ are isomorphic if and only if $\lmd'=\lmd$ and
 $G'=\xi G$ for some nonzero $\xi\inc$.
 Moreover, any Lie algebra isomorphism $\pi:\g(G,\lmd)\longrightarrow\g(\xi G,\lmd)$,
 aside from an inner automorphism of $\g(\xi G,\lmd)$, has the form
 \begin{equation}\label{eq4.01}
 \begin{aligned}
  &\pi(L_a)=\xi^{-1}\chi(a)L_{\xi a}'+\chi(a)I_{\xi a}'\left(
            l_0\dt_{\lmd,0}+l_1a+l_2a^2\dt_{\lmd,-1}+l_3a^3\dt_{\lmd,-2}+f(a)\dt_{\lmd,1}\right);\\
  &\pi(I_a)=l\chi(a)I_{\xi a}',
 \end{aligned}
 \end{equation}
 where $l\inc^\times,l_0,l_1,l_2,l_3\in\C$, $\chi\in G^*$ is a character of $G$,
 and $f\in\text{Hom}(G,\C)$.
\end{thm}
\pf{
Denote $\g=\g(G,\lmd), \g'=\g(G',\lmd')$,
and we will use an extra dash to denote elements in $\g'$ in the following.
If $\lmd'=\lmd,\ G'=\xi G$ for some nonzero $\xi\inc$,
it is easy to check that the linear map defined by
$$L_a\mapsto \xi^{-1}L_{\xi a}';\ \ \ I_a\mapsto \xi^{-1}I_{\xi a}'$$
is a Lie algebra isomorphism from $\g$ to $\g'$.

On the other hand, suppose $\g\cong\g'$ and
let $\pi:\g\longrightarrow\g'$ be a Lie algebra isomorphism.

{\bf Claim 1}: There exists some $\xi\in\C^\times$
and an inner automorphism $\theta$ of $\g'$ such that $G'=\xi G$ and
$\theta\pi(\g_a)=\g_{\xi a}'$ for any $a\ing$.

Notice that the set of locally finite elements in $\g$ is $\spanc{L_0,I_a\mid a\ing}$.
Since $\pi$ maps a locally finite element in $\g$ to a locally finite element in $\g'$,
we may assume that
$$\pi(L_0)=\xi^{-1}L_0'+\sum_{b\in\setgl 0}\gm_bI_b',$$
for some $\xi\in\C^\times$, finite subset $\setgl0$ of $G'$ and $\gm_b\inc$.
Define
$$\eta=\prod_{0\neq b\in\setgl0}\exp{\left\{-\frac{\xi\gm_b}b\text{ad}I_b'\right\}},$$
which is an inner automorphism of $\g'$.
Clearly, $\eta\pi(L_0)=\xi^{-1}L_0'+\gm_0I_0'.$
So replacing $\pi$ by $\eta^{-1}\pi$ we may assume
$$\pi(L_0)=\xi^{-1}L_0'+\gm_0I_0'.$$
Since $\spanc{I_a\mid a\ing}$ is the unique maximal abelian ideal of $\g$, it forces
$$\pi(I_a)=\sum_{b\in\setgi a}\nu_a(b)I_b'$$
for some finite subset $\setgi a$ of $G'$ and some function $\nu_a:G'\longrightarrow\C$.

For any $a\neq 0$, from
$$\pi[L_0,I_a]=a\sum_{b\in\setgi a}\nu_a(b)I_b'=[\pi(L_0),\pi(I_a)]
  =\xi^{-1}\sum_{b\in\setgi a}b\nu_a(b)I_b'$$
we see that $\setgi a=\{\xi a\}\subset G'$ for any $a\neq0$.
This implies $G'=\xi G$ and
$$\pi(I_a)=\nu_a(\xi a)I_{\xi a}'\text{ for } a\neq0.$$
Since $\pi[L_a,I_b]=(b-\lmd a)\pi(I_{a+b})=(b-\lmd a)\nu_{a+b}(\xi(a+b))I_{\xi(a+b)}'
        =[\pi(L_a),\nu_b(\xi b)I_{\xi b}']$
for $a,b\neq0,a+b\neq0$,
it follows that $\pi(L_a),a\neq0,$ must have the form
$$\pi(L_a)=\rho(a)L_{\xi a}'+\sum_{c\in\setgl a}\mu(c)I_c'$$
for some finite subset $\setgl a$ of $G'$, functions $\rho:G\longrightarrow\C^\times$
and $\mu:G'\longrightarrow\C$. Then
\begin{equation}\label{eq4.02}
 (b-\lmd a)\nu_{a+b}(\xi(a+b))I_{\xi(a+b)}'=[\rho(a)L_{\xi a}',\nu_b(\xi b)I_{\xi b}']
  =\xi\rho(a)\nu_b(\xi b)(b-\lmd'a)I_{\xi(a+b)}'.
\end{equation}
Moreover, from
$$\pi[L_a,I_0]=-\lmd a\pi(I_a)=-\lmd a\nu_a(\xi a)I_{\xi a}'
  =[\rho(a)L_{\xi a}',\sum\limits_{b\in\setgi0}\nu_0(b)I_b']$$
we see that $\setgi 0=\{0\}\subset G'$, and hence $\pi(I_0)=\nu_0(0)I_0'$.
Define
$$\nu(\xi b)=\nu_b(\xi b)\text{ for any }b\ing.$$
We get a function $\nu:G'\longrightarrow\C^\times$ such that
$$\pi(I_a)=\nu(\xi a)I_{\xi a}'\in\g'_{\xi a}\ \ \ \text{ for any }a\ing.$$
Then (\ref{eq4.02}) turns to
\begin{equation}\label{eq4.03}
 (b-\lmd a)\nu(\xi(a+b))=\xi\rho(a)\nu(\xi b)(b-\lmd'a).
\end{equation}
Consider $\pi[L_0,L_a]$ for $a\neq 0$, we get
$$\sum_{b\in\setgl a}\mu(b)(a-\xi^{-1}b)I_b'-\gm_0\lmd'\xi a\rho(a)I_{\xi a}'=0,$$
which implies $\setgl a=\{\xi a\}$ and $\lmd'\gm_0=0$.
Write $\mu(0)=\gm_0$ and $\rho(0)=\xi^{-1}$. Then we have
$$\pi(L_a)=\rho(a)L_{\xi a}'+\mu(\xi a)I_{\xi a}'\in\g'_{\xi a}\ \ \ \text{ for any }a\ing.$$
This proves Claim 1.

In the following we determine the functions $\rho,\mu,\nu$.\\
{\bf Claim 2}: $\chi=\xi\rho$ is a character of $G$,
$\lmd'=\lmd$ and $\nu(\xi a)=\nu(0)\chi(a)$ for any $a\ing$.
Since
$$\begin{aligned}
 &\pi[L_a,L_b]=(b-a)\pi(L_{a+b})=(b-a)\rho(a+b)L_{\xi(a+b)}'+(b-a)\mu(\xi(a+b))I_{\xi(a+b)}'\\
 &=[\rho(a)L_{\xi a}'+\mu(\xi a)I_{\xi a}',\rho(b)L_{\xi b}'+\mu(\xi b)I_{\xi b}']\\
 &=\xi(b-a)\rho(a)\rho(b)L_{\xi(a+b)}'+\rho(a)\mu(\xi b)\xi(b-\lmd'a)I_{\xi(a+b)}'
   -\rho(b)\mu(\xi a)\xi(a-\lmd'b)I_{\xi(a+b)}',
\end{aligned}$$
we get
\begin{align}
 &(b-a)\left(\rho(a+b)-\xi\rho(a)\rho(b)\right)=0,\label{eq4.04}\\
 &(b-a)\mu(\xi(a+b))=\xi\rho(a)\mu(\xi b)(b-\lmd'a)-\xi\rho(b)\mu(\xi a)(a-\lmd'b).\label{eq4.05}
\end{align}
From (\ref{eq4.04}) we get that if $a\neq b$ then
\begin{equation}\label{eq4.06}
 \rho(a+b)=\xi\rho(a)\rho(b).
\end{equation}
Let $b=0$ in (\ref{eq4.06}), we have $\rho(0)=\xi^{-1}$.
Choose $b\neq 0,\pm a$. Using (\ref{eq4.06}) we get
$$\rho(2a)=\xi\rho(a+b)\rho(a-b)=\xi^3\rho(a)^2\rho(b)\rho(-b)
  =\xi^2\rho(a)^2\rho(0)=\xi \rho(a)^2\ \text{  for any }a\ing.$$
This proves that (\ref{eq4.06}) stands for any $a,b\ing$.
So $\chi=\xi\rho$ is a character of $G$.

Let $b=0$ in (\ref{eq4.03}) we get
\begin{equation}\label{eq4.07}
 \lmd\nu(\xi a)=\lmd'\xi\nu(0)\rho(a)\text{ for any }a\neq0,
\end{equation}
which implies that
$$\lmd'=0\text{ if and only if }\lmd=0.$$
If $\lmd\neq0$, then (\ref{eq4.07}) says
$\nu(\xi a)=\frac{\lmd'\xi}\lmd\nu(0)\rho(a)$ for any $a\neq0$.
Putting it into (\ref{eq4.03}) and using (\ref{eq4.06}) we see $\lmd'=\lmd$.
Hence $\nu(\xi a)=\xi\nu(0)\rho(a)=\nu(0)\chi(a)$ if $\lmd\neq0$.

If $\lmd=\lmd'=0$, let $a=-b\neq0$ in (\ref{eq4.03}) and we have
$$\nu(0)=\xi\rho(-b)\nu(\xi b).$$
Then $\nu(0)\chi(b)=\xi^2\rho(b)\rho(-b)\nu(\xi b)=\xi\rho(0)\nu(\xi b)=\nu(\xi b)$.
This proves Claim 2.

At last we determine the function $\mu$.
Set $\vf(a)=\frac{\mu(\xi a)}{\chi(a)}$.
Divide $\chi(a+b)$ to (\ref{eq4.05}) and we obtain
$$(b-a)\vf(a+b)+(a-\lmd b)\vf(a)=(b-\lmd a)\vf(b),$$
which is the same equation as (\ref{eq3.04}).
So from the computation in Section 3 we have
$$\vf(a)=l_0\dt_{\lmd,0}+l_1a+l_2a^2\dt_{\lmd,-1}+l_3a^3\dt_{\lmd,-2}+f(a)\dt_{\lmd,1},$$
where $l_0,l_1,l_2,l_3\inc$ and $f\in\text{Hom}(G,\C)$.
Therefore
$$\mu(\xi a)=\chi(a)\left(
  l_0\dt_{\lmd,0}+l_1a+l_2a^2\dt_{\lmd,-1}+l_3a^3\dt_{\lmd,-2}+f(a)\dt_{\lmd,1}\right).$$
This proves the theorem.
}

Now using Theorem \ref{thm4.1} we may determine the automorphism group of $\g$.
Denote $\mathcal E=\{\e\inc^\times\mid \e G=G\}$, which is a subgroup of $\C^\times$.
Let $\xi\in\mathcal E, \chi\ing^*,f\in\text{Hom}(G,\C),l\inc^\times,l_0,l_1,l_2,l_3\inc$,
denote by $\theta_\lmd(\xi, \chi,f,l,l_0,l_1,l_2,l_3)$ the linear map given by
\begin{align*}
  &L_a\mapsto\xi^{-1}\chi(a)L_{\xi a}+\chi(a)I_{\xi a}\left(
            l_0\dt_{\lmd,0}+l_1a+l_2a^2\dt_{\lmd,-1}+l_3a^3\dt_{\lmd,-2}+f(a)\dt_{\lmd,1}\right);\\
  &I_a\mapsto l\chi(a)I_{\xi a}.
\end{align*}
Here we shall point out that if $\lmd=1$ we may assume $l_1=0$
since the map $a\mapsto f(a)+l_1a$ still lies in $\text{Hom}(G,\C)$.
Then by Theorem \ref{thm4.1} we have
\begin{thm}
 The outer automorphism group of $\g$ is
 $$\text{Out}\g=\{\theta_\lmd(\xi, \chi,f,l,l_0,l_1,l_2,l_3)\mid
 \xi\in\mathcal E, \chi\ing^*,f\in\text{Hom}(G,\C),l\inc^\times,l_0,l_1,l_2,l_3\inc\},$$
 and $\text{Aut}\g=\text{Out}\g\ltimes\text{Inn}\g$, where
 $\text{Inn}\g$ is the inner automorphism group of $\g$, generated by $\{\exp\ad I_a\mid a\ing\}$.
\end{thm}

Notice that
\begin{align*}
 &\theta_\lmd(\xi',\chi',f',l',l_0',l_1',l_2',l_3')\cdot\theta_\lmd(\xi, \chi,f,l,l_0,l_1,l_2,l_3)\\
 &=\theta_\lmd\left(\xi'\xi,(\chi'\cdot\xi)\chi,f'+l'f,l'l,
                \xi^{-1}l_0'+l'l_0,l_1'+l'l_1,\xi l_2'+l'l_2,\xi^2l_3'+l'l_3\right)
\end{align*}
and
\begin{align*}
 \theta_\lmd(\xi,&\chi,f,l,l_0,l_1,l_2,l_3)^{-1}=\\
        &\theta_\lmd\left(\xi^{-1},\chi^{-1}\cdot\xi^{-1},-l^{-1}f,l^{-1},
          -l^{-1}\xi l_0,-l^{-1}l_1,-l^{-1}\xi^{-1}l_2,-l^{-1}\xi^{-2}l_3\right).
\end{align*}
Clearly, the map $G^*\longrightarrow\text{Out}\g$ defined by
$\chi\mapsto \theta_\lmd(1,\chi,0,1,0,0,0,0)$
is a group monomorphism and denote by $N$ its image.
Set $S=\{\theta_\lmd(1,\mathbbm{1},0,l,0,0,0,0)\mid l\inc^\times\}$, and
$$K=\{\theta_\lmd(1,\mathbbm{1},f,1,l_0,l_1,l_2,l_3)\mid
        l_0,l_1,l_2,l_3\inc,f\in\text{Hom}(G,\C)\}.$$
Here $\mathbbm 1$ denotes the identity in $G^*$.
Clearly, $N,S,K$ are normal subgroups of $\text{Out}\g$, and
$$N\cong G^*,\
S\cong\C^\times, \text{ and }
K\cong\begin{cases}
         \text{Hom}(G,\C)\cong\C^n &\text{ if }\lmd=1;\\
         \C^2 &\text{ if }\lmd=0,-1,-2;\\
         \C   &\text{ otherwise}.
      \end{cases}$$
Moreover, let $T=\{\theta_\lmd(\xi,\mathbbm{1},0,1,0,0,0,0)\mid \xi\in\mathcal E\}$,
which is a subgroup of $\text{Out}\g$ and isomorphic to $\mathcal E$.
Then we have a projection $\text{Out}\g\longrightarrow T$,
whose kernel is the normal subgroup $NSK$.
So we get
\begin{thm}
 $\text{Out}\g=T\ltimes(NSK)\cong\mathcal E\ltimes\big(G^*\C^\times K\big).$
\end{thm}

Using a result from \cite{P} one may extend automorphisms of $\g$
to automorphisms of $\overline\g$.

\begin{prop}[\cite{P}]
 Let $\G$ be a perfect Lie algebra and $\overline\G$ the universal central extension of $\G$.
 Then every automorphism $\theta$ of $\g$ admits a unique extension to an automorphism $\overline\theta$
 of $\overline\G$, and the map $\theta\mapsto\overline\theta$ is a group monomorphism.
 Moreover, if $\G$ is centerless,
 then the map $\theta\mapsto\overline\theta$ is an isomorphism.
\end{prop}

Explicitly we have the form of extended automorphisms.
\begin{thm}
 Suppose $\lmd\neq-1$.
 Let $\xi\in\mathcal E, \chi\ing^*,f\in\text{Hom}(G,\C),l\inc^\times,l_0,l_1,l_3\inc$,
 then the unique automorphism $\ov\theta=\ov{\theta_\lmd(\xi,\chi,f,l,l_0,l_1,0,l_3)}$ of $\ov\g$
 obtained by extension from $\theta=\theta_\lmd(\xi,\chi,f,l,l_0,l_1,0,l_3)$
 is such that
  \begin{align}
    &\ov\theta(L_a)=\xi^{-1}\chi(a)L_{\xi a}  \label{eq4.08}
       +\chi(a)I_{\xi a}\left(l_0\dt_{\lmd,0}+l_1a+l_3a^3\dt_{\lmd,-2}+f(a)\dt_{\lmd,1}\right)\\
     &\hspace{2cm}+\dt_{a,0}\left(\frac{\xi^{-1}-\xi}{24}C_L+(l_1\xi-l_0)\cli0\dt_{\lmd,0}
         +\frac12\xi(l_1^2-l_0^2)C_I\dt_{\lmd,0}\right); \notag  \\
      &\ov\theta(I_a)=l\chi(a)I_{\xi a}
        +l\xi\dt_{a,0}\left(\frac{\xi^{-1}-\xi}{24}\cli1\dt_{\lmd,1}
                 +\left((1-\xi^{-1})\cli0+(l_1-l_0)C_I\right)\dt_{\lmd,0}\right)\label{eq4.09}\\
     &\hspace{2cm}+l(\xi\e_1)^{-1}\dt_{a,0}\sum_{i=2}^n(\xi\e_1)_i\cli i\dt_{\lmd,-2}; \notag \\
     &\ov\theta(C_L)=\xi C_L-12l_1\xi\left(2\cli0+l_1C_I\right)\dt_{\lmd,0};\label{eq4.10}\\
     &\ov\theta(\cli0)=l\xi\left(\cli0+l_1C_I\right);
          \hspace{1cm}\ov\theta(\cli1)=l\xi^2 C_{LI}^{(1)};\hspace{1cm}  \notag \\
     &\hspace{1cm}\ov\theta(\cli i)=l(\xi\e_1)^{-1}\sum_{j=2}^n
          \left\{\e_1(\xi\e_i)_j-\e_i(\xi\e_1)_j\right\}\cli j,\ i\geq2;\label{eq4.11}\\
     &\ov\theta(C_I)=l^2\xi C_I.\label{eq4.12}
  \end{align}
 Here for $a\in G$ we use $a_i$ to denote the coefficients of $a$ with respect to the basis
 $\e_1,\dots,\e_n$.
\end{thm}
\pf{
Recall $\mathfrak c=\spanc{C_L,C_I\dt_{\lmd,0},\cli 0\dt_{\lmd,0},\cli 1\dt_{\lmd,1},
               \cli i\dt_{\lmd,-2}\mid 2\leq i\leq n}$.
For later convenience, we denote
\begin{align*}
 &\tau(a)=l_0\dt_{\lmd,0}+l_1a+l_3a^3\dt_{\lmd,-2}+f(a)\dt_{\lmd,1},\\
 &\nu_C(a)=C_{LI}^{(0)}(a^2+a)\dt_{\lmd,0}+\frac1{12}(a^3-a)C_{LI}^{(1)}\dt_{\lmd,1}
   +\sum\limits_{i=2}^na_iC_{LI}^{(i)}\dt_{\lmd,-2}\in\mathfrak c.
\end{align*}
Clearly, $\ov\theta(L_a)$ and $\ov\theta(I_a)$ have the form
\begin{align*}
 \ov\theta(L_a)=&\xi^{-1}\chi(a)L_{\xi a}+\chi(a)\tau(a)I_{\xi a}+K_L(a);\\
 \ov\theta(I_a)=&l\chi(a)I_{\xi a}+K_I(a),
\end{align*}
for some maps $K_L,K_I:G\longrightarrow\mathfrak c$.
Expanding
\begin{align*}
   \ov\theta([L_a,L_b])=[\xi^{-1}\chi(a)L_{\xi a}+\chi(a)\tau(a)I_{\xi a},
            \xi^{-1}\chi(b)L_{\xi b}+\chi(b)\tau(b)I_{\xi b}]
\end{align*}
we see that $K_L(a)=0$ if $a\neq0$, and
$$-2K_L(0)+\frac1{12}(a^2-1)\ov\theta(C_L)=\frac{\xi a^2-\xi^{-1}}{12}C_L
        +2(l_0-l_1\xi a^2)\cli0\dt_{\lmd,0}+(l_0^2-l_1^2a^2)\xi C_I\dt_{\lmd,0},$$
which implies
$$K_L(0)=\frac{\xi^{-1}-\xi}{24}C_L+(l_1\xi-l_0)\cli0\dt_{\lmd,0}
        +\frac12\xi(l_1^2-l_0^2)C_I\dt_{\lmd,0}$$
and
$$\ov\theta(C_L)=\xi C_L-12l_1\xi\left(2\cli0+l_1C_I\right)\dt_{\lmd,0}.$$
This proves (\ref{eq4.08}) and (\ref{eq4.10}).

Expanding $\ov\theta([L_a,I_b])=[\ov\theta(L_a),\ov\theta(I_b)]$ we get that
$K_I(a)=0$ for $a\neq0$ and
\begin{equation}\label{eq4.13}
 \begin{split}
  &-(1+\lmd)aK_I(0)+(a^2+a)\ov\theta(\cli0)\dt_{\lmd,0}
     +\frac{a^3-a}{12}\ov\theta(\cli1)\dt_{\lmd,1}
     +\sum_{i=2}^na_i\ov\theta(\cli i)\dt_{\lmd,-2}\\
  &\ =l\xi^{-1}\left(\left((\xi a)^2+\xi a\right)\cli0\dt_{\lmd,0}
         +\frac{(\xi a)^3-\xi a}{12}\cli1\dt_{\lmd,1}
       +\xi\sum_{i=2}^na_i\cli i\dt_{\lmd,-2}\right)
       +l\tau(a)\xi aC_I\dt_{\lmd,0}.
 \end{split}
\end{equation}
If $\lmd=0$, let $a=-1$ in (\ref{eq4.13}), then we get
$$K_I(0)\dt_{\lmd,0}=l\xi\left((1-\xi^{-1})\cli0+(l_1-l_0)C_I\right)\dt_{\lmd,0}\ \text{ and }\
  \ov\theta(\cli0)=l\xi\left(\cli0+l_1C_I\right).$$
If $\lmd=1$, let $a=1$ in (\ref{eq4.13}), then we get
$$K_I(0)\dt_{\lmd,1}=-\frac{l(\xi^2-1)}{24}\cli1\dt_{\lmd,1}\ \text{ and }\
  \ov\theta(\cli1)=l\xi^2\cli1.$$
If $\lmd=-2$, let $0\neq a\in\Z\e_1$ in (\ref{eq4.13}), then we get
\begin{align*}
 &K_I(0)\dt_{\lmd,-2}=l(\xi\e_1)^{-1}\sum_{i=2}^n(\xi\e_1)_i\cli i\dt_{\lmd,-2};\\
 &\ov\theta(\cli i)=l(\xi\e_1)^{-1}\sum_{j=2}^n
          \left(\e_1(\xi\e_i)_j-\e_i(\xi\e_1)_j\right)\cli j.
\end{align*}
This proves (\ref{eq4.09}) and (\ref{eq4.11}).
Similarly, (\ref{eq4.12}) follows from
$\ov\theta([I_a,I_b])=[\ov\theta(I_a),\ov\theta(I_b)]$.
}

\end{document}